\documentclass[12pt]{article}

\usepackage{andieng}
\usepackage{amsmath}
\usepackage{dsfont}
\usepackage{amssymb}
\usepackage{amsthm}
\def\R{\mathds{R}}

\def\H{\mathds{H}}
\def\N{\mathds{N}}
\def\Z{\mathds{Z}}
\def\C{\mathds{C}}

\def\Kr#1{      \left( #1 \right) }
\def\st{\,|\,}

\usepackage{graphicx}
\usepackage{hyperref}

\newtheorem{lemma}{Lemma}[section]
\newtheorem{corollary}[lemma]{Corollary}
\newtheorem{theorem}[lemma]{Theorem}
\newtheorem{proposition}[lemma]{Proposition}

\newtheorem*{mainthm}{Main Theorem}
\newtheorem*{COR}{Corollary}
\newtheorem*{proposition*}{Proposition}

\theoremstyle{definition}
\newtheorem{definition}[lemma]{Definition}
\newtheorem{remark}[lemma]{Remark}

\TopRand{0cm}
\def\slope{{\rm slope}}

\def\EB{Euclidean building}
\def\HS{Hadamard space}
\def\C{\mathcal{G}}

\def\comment#1{}

\def\color[#1]#2{}

\begin{document}
\title{Polygons with prescribed Gauss map in Hadamard spaces and Euclidean
buildings}
\author{Andreas Balser}
\date{December 2, 2004}

\maketitle
\begin{abstract}
\noindent
We show that given a stable weighted configuration on the asymptotic
boundary of a 
locally compact Hadamard space, there is a polygon with Gauss
map prescribed by the given weighted configuration if the
configuration is stable. Moreover, the same result holds for
semistable configurations on arbitrary \EB s.
\footnote{MSC2000: 53C20\\
\emph{keywords:} \EB s, \HS s, polygons}
\end{abstract}

\noindent
In the first section, we recall some background material on Hadamard
spaces and \EB s, and we introduce the concepts needed to state and
prove our  
Theorems. In particular, we define stability for weighted
configurations on the boundary at infinity of a Hadamard space.

In the second section, we introduce ultralimits and the special cases
ultraproducts and asymptotic tubes which we use in our proofs.

In the third section, we prove our results:

\begin{mainthm}
Let $X$ be a \EB\ and $c$ a
semistable weighted configuration on its boundary at infinity, or let
$X$ be a locally 
compact \HS\ and $c$ a stable weighted configuration on its boundary
at infinity.
Then the  
associated weak contraction~$\Phi_{c}$ has a fixed point.
In particular, there exists a polygon
$p$ in $X$ such that $c$ is a Gauss map for $p$.
\end{mainthm}

For a slightly more general statement in the case of a \HS , see
Corollary~\ref{hadamardCorollary}.

As an immediate consequence, we can formulate the following
classification of configurations which can occur as Gauss maps on
Euclidean buildings and symmetric spaces:

\begin{COR}
Let $X$ be a symmetric space of non-compact type or a \EB, and let $c$ be a
weighted configuration on its boundary at infinity. 

Then there exists a polygon having this configuration as a Gauss map
if and only if 
the configuration is semistable in the building case and nice
semistable in the case of a symmetric space.
\end{COR}

Necessity of
semistability, as well as the Theorem and the Corollary  in the case
where~$X$ is a symmetric space 
or a locally compact \EB \ were shown in \cite{KLM1},
\cite{kapovichLeebMillson}.

We extend their ideas by suitable use of ultralimits.

Along the way, we discuss how rays project to subspaces, and obtain
the following result of independent interest:

\begin{proposition*}
Let $X'$ be a Hadamard space and $X\subset X'$ a closed convex
subset.  Let
$o\in X$, $\rho :=\overline{o\eta}$ a ray in $X$ (with $\eta \in
\partial_{\infty}X'$), and $\pi:X'\rightarrow X$ be the 
nearest point
projection. 

If $\angle (\eta ,\partial_{\infty}X)<\frac{\pi}{2}$,
then the segments $\overline{o\,\,(\pi \circ \rho (t))}$
converge to the 
ray $\overline{o\xi}$ (in the cone topology), where $\xi\in
\partial_{\infty}X$ is the unique point with $\angle (\eta ,\xi
)=\angle (\eta,\partial_{\infty}X)$.

If the projection
 of the ray
$\overline{o\eta}$ to $X$ is bounded, 
then there exists a point $p\in X$ s.t.\ $\pi \circ
\overline{p\eta}(t)=p$ for all $t>0$.
\end{proposition*}

In the last section, we discuss relations of the above results to
algebra.

The author wishes to thank Vitali Kapovich and Viktor Schroeder for
the suggestion to think about the generalization to \HS s, and
Robert Kremser, Bernhard Leeb, Alexander Lytchak and John Millson for
discussions 
about the topic. 

Further, I would like to point out that upon finishing this article,
I found out that Misha Kapovich had independently had the same idea for 
proving the main theorem for buildings. 

\section{Introduction}
\subsection{\HS s}
We will use the language of non-positively curved metric spaces, as
developed in \cite{ballmann}. 

Throughout, let $X$ be a Hadamard space, unless otherwise stated.

Recall that $X$ has a \emph{boundary at infinity}
$\partial_{\infty}X$, which is given by equivalence classes of rays,
where two (unit-speed) rays are equivalent if their distance is
bounded. 

In particular, we will use Busemann functions $b_{\eta}$ associated to
an asymptotic boundary point $\eta \in \partial_{\infty}X$.
A Busemann function measures (relative) distance from a point at
infinity, and is determined up to an additive constant
only. Busemann functions are convex (along any geodesic) and
1-Lipschitz. 

Geodesics, rays, and geodesic segments are always assumed to be
parametrized by unit speed (i.e.\ they are isometric embeddings). 

For a line $l$ in $X$, there is the space $P_{l}$ of parallel
lines. $P_{l}$ splits as a product $P_{l}\cong l\times CS(l)$, where
$CS(l)$ is a \HS \ again. 

For points $x,\xi$ with $x\in X$, $\xi \in X\cup \partial_{\infty}X$,
and $t\geq 0$ (if $\xi \in X$, let $t\leq d(x,\xi )$ ), 
we let $\overline{x\xi}(t)$ denote the point on the segment/ray
$\overline{x\xi}$ at distance $t$ from $x$. When we denote a ray by
$\overline{o\eta}$, we order the points such that $o\in X$ and $\eta
\in \partial_{\infty}X$. 

\begin{definition}
For $\xi \in \partial_{\infty}X$ and $t\geq 0$, we define the map
$\phi_{\xi ,t}:X\rightarrow X$ defined by 
$\phi_{\xi,t}(x):=\overline{x\xi}(t)$. 
Observe that $\phi_{\xi ,t}$ is a 1-Lipschitz map by convexity of a
non-positively curved metric.
\end{definition}

Let $o\in X$ be a point in a \HS , and let $\eta ,\xi \in
\partial_{\infty}X$. Let $c,c'$ be the rays
$\overline{o\eta},\overline{o\xi}$. For points $c(t),c'(t')$, one can
consider the \emph{Euclidean comparison triangle} corresponding to the
points $o,c(t),c'(t')$, i.e.\ the Euclidean triangle with side-lengths
$d(o,c(t)), d(c(t),c'(t')), d(c'(t'),o)$ (which is well-defined up to
isometries of the Euclidean plane). The \emph{comparison angle}
between $c(t)$ and $c'(t')$ at $o$ is the angle of the comparison
triangle at the point corresponding to $o$. It is denoted by
$\tilde{\angle}_{o}(c(t),c'(t'))$. 

We have the following monotonicity property: 
\[
0<t\leq s\text{ and }0<t'\leq s'\text{ implies }
\tilde{\angle}_{o}(c(t),c'(t'))\leq
\tilde{\angle}_{o}(c(s),c'(s')).
\]

From this, one can deduce a notion of angle between geodesic segments
and rays: 
\[
\angle_{o}(\eta,\xi  ):=\lim_{t,t'\rightarrow
0}\tilde{\angle}_{o}(c(t),c'(t'))\in [0,\pi ],
\]
and an ``angle at infinity'', the \emph{Tits angle} between boundary
points
\[
\angle (\eta,\xi  ):=\angle_{Tits}(\eta ,\xi ):=\lim_{t,t'\rightarrow
\infty}\tilde{\angle}_{o}(c(t),c'(t'))	\in [0,\pi ].
\]
It is easy to see that the Tits angle between $\eta ,\xi$ does not
depend on the chosen 
basepoint~$o$. 
The length metric induced on $\partial_{\infty}X$ by $\angle$ is
called \emph{Tits distance} $Td$, and makes $\partial_{\infty}X$ a
CAT(1) space.
If the Tits angle (between $\eta ,\xi$) is less than $\pi$, there is a 
unique geodesic $\overline{\eta \xi}\subset \partial_{\infty}X$
connecting them. 

Similarly, the space of directions $S_{o}$, i.e.\ the completion of
the space 
of starting directions of geodesic segments initiating in
$o$ (modulo the equivalence of directions enclosing a zero angle), can
be regarded as a CAT(1) space.

\medbreak
We state Lemma \cite[II.8.3]{bridsonHaefliger}, since it will be of
fundamental importance in the proof of Lemma \ref{commutativity}. 
It says that given one geodesic ray $\overline{o\eta}$ and another
point $y\in X$, the ray $\overline{y\eta}$ can be approximated by
segments $\overline{y\rho (t)}$ for $t$ large enough, and the
approximation can be controlled
independently from the \HS ~$X$.

\begin{lemma}\label{rayLemma}
Given $\varepsilon >0, m>0$ and $c>0$, there is a constant
$K=K(\varepsilon ,m,c)>0$ such that: Let $\rho$ be a ray
$\overline{o\eta}$ in a
Hadamard-space $X$. If $y\in X$ satisfies $d(y,o)\leq m$, then
we have 
\begin{align}
d(\overline{y\eta}(c), \overline{y\rho (K)}(c))<\varepsilon.\tag*{\qed}
\end{align}
\end{lemma}

\subsection{\EB s}
We will also need some Euclidean building geometry. For an
introduction, we refer to \cite[sect.~4]{kleinerLeeb}.
A brief introduction of the notation we use can be found in
\cite[sect.~2.4]{kapovichLeebMillson}.
Note that in particular, a \EB \ is a Hadamard space.

The boundary at infinity of a \EB \ $X$ of rank $n$ is a spherical
building of dimension $n-1$; 
we refer to \cite[sect.~3]{kleinerLeeb} for an introduction.

We will use that a spherical building is a spherical simplicial
complex, where all the simplices are isometric to a spherical polytope
$\Delta$ (in particular, $\Delta$ tesselates $S^{n-1}$), which is the
\emph{spherical Weyl chamber} of the building. Apartments
(i.e.\ isometrically embedded copies $S^{n-1}$) intersect in (unions
of) Weyl
chambers.

We prove some elementary lemmas which we will use later:

\begin{lemma}\label{parallels}
Let $X$ be a \EB , $l$ a line in $X$ with $l(\infty )=\eta \in
\partial_{\infty}X$, and $c$ a ray asymptotic to $\eta$. Then $c$
eventually coincides with a line parallel to $l$.
\end{lemma}
\begin{proof}
Pick an apartment $A'\supset c$, and an apartment $A$ containing
$\eta^{-}:=l(-\infty )$ in its boundary, which has the property that
$\partial A=\partial A'$ near $\eta$ (i.e.: let $S\subset
\partial_{\infty}A'$ be the subset of $\partial_{\infty}A'$ consisting
of the union of Weyl chambers containing $\eta$, and let $A$ be an apartment
containing $S$ and $\eta^{-}$ in its boundary).  

We want to show that $c(t)\in A$
for large $t$, which finishes the proof.

We may assume that  $ \eta $ is singular, since otherwise $c(t)\in A$
for large $t$ by \cite[L.~4.6.3]{kleinerLeeb}.

Pick regular points $\xi_{i}\in S$
such that $\eta$ is the midpoint of $\overline{\xi_{1}\xi_{2}}$ (and
$\angle_{Tits} (\xi_{1},\xi_{2})<\pi$).

Let $c_{i}$ be the ray $\overline{c(0)\xi_{i}}\subset A'$. For some $t_{0}$, both
$c_{i}(t_{0})\in A\cap A'$ (again by
\cite[L.~4.6.3]{kleinerLeeb}). Then the midpoint of  
$\overline{c_{1}(t_{0})c_{2}(t_{0})}$ is also in $A\cap A'$; this
midpoint is $c(T)$ for some $T$ (since $c_{1},c_{2}$ span a flat
sector in $A'$), implying
that $c(t)$ is in $A\cap A'$ for $t\geq T$. 
\end{proof}

Observe that this shows in particular that the space of strong
asymptote classes of rays asymptotic to $\eta$ is isometric to
$CS(l)$
(see \cite{karpelevic}, \cite[sect.\ 2.1.3]{leebHabil}).

\begin{lemma}\label{raysApartments}
Consider a ray $\rho =\overline{o\xi}$ and a segment $\overline{op'}$ in
a \EB \ $X$. Then there is an apartment containing $\rho$ and an
initial part of $\overline{op'}$.  
\end{lemma}
\begin{proof}
The claim is clear if $\rho$ and $\overline{op'}$ initially coincide
or their initial directions are antipodal.
So we assume that they do not.
By \cite[L.\ 4.1.2]{kleinerLeeb}, there is a point $p\in
\overline{op'}$, such that the triangle $D:=\Delta (o,p,\xi)$ is flat
(i.e.\ a flat 
half-strip). If $X$ is discrete, the claim follows from
\cite[Prop.\ 1.3, Rem.\ 1.4]{buildingLike}. We give a direct argument
for our special situation here:

We show that $D$ is contained in a half-plane: Let $H$ be a flat
half-strip containing~$D$ with $\partial H\supset \overline{op}$;
assume that $H$ cannot be enlarged under these conditions, and is not
a half-plane. 
Since $X$ is complete, we see that $H$ is closed, i.e.\ of the form
$\Delta (p_{1},p_{2},\xi)$. Now $S_{p_{1}}H$ is a geodesic segment,
which can be prolonged to a geodesic of length $\pi$ in the spherical building
$S_{p_{1}}X$. By \cite[L.\ 4.1.2]{kleinerLeeb}, this yields a
direction in which we can glue another flat half-strip to $H$, so $H$
was not maximal.

Thus, $D$ is contained in a half-plane, and this half-plane is
contained in a plane by \cite[L.\ 5.2]{leebHabil}. Finally, every
plane in $X$ is contained in an apartment by \cite[Cor.~5.4]{leebHabil}.
\end{proof}

\subsection{Weighted configurations at infinity}

In this subsection, we recall some  notions from \cite{KLM1} and
\cite{kapovichLeebMillson} needed to discuss the 
relationship of configurations on $\partial_\infty X$ and polygons in $X$.

\begin{definition}
Let $X$ be a \HS . A \emph{weighted configuration $c$ on $\partial_\infty X$} is
an $n$-tuple of points $(\xi _{1},\dotsc ,\xi _{n})$ in $\partial_\infty X$ together
with a weight function $m:\{1,\dotsc ,n \}\rightarrow \R_{> 0}$.
\end{definition}

There is a \emph{weighted Busemann function} associated to a weighted
configuration $c$. It is given by
\[
b_{c}:=\sum_{i=1}^{n}m_{i}b_{\xi_{i}};
\]
weighted Busemann functions are convex, asymptotically linear,
Lipschitz-continuous, and 
well-defined up to an additive 
constant.
As for any convex, asymptotically linear Lipschitz-function on a \HS , we can
associate a  
function $\slope_{b_{c}}: \partial_\infty X\rightarrow \R$ to a weighted
Busemann function, which is given by assigning the asymptotic slope of
$b_{c}$ on a ray $\overline{o\xi}$ to the point $\xi$. Since two rays
asymptotic to the same boundary point have bounded distance and
$b_{c}$ is Lipschitz, 
 the slope does not depend on the choice of $o$, so
$\slope_{b_{c}}$ is well-defined (see also \cite[sect.~3]{KLM1}).

We have 
\[
\slope_{c}(\xi ):=\slope_{b_{c}}(\xi)=
-\sum_{i=1}^{n}m_{i}\cos\angle (\xi_{i},\xi ).
\]

The configuration $c$ is called \emph{semistable} if
$\slope_{c}\geq 0$, and
it is called \emph{stable} if $\slope_{c}>0$. 

Observe that (semi-)stability is defined purely in terms of the
Tits-geometry of $\partial_{\infty}X$, without reference to $X$ itself.

\medbreak

Now we discuss the relation between polygons and weighted
configurations: 

Consider a polygon $p$ in $X$, which is 
determined by an $n$-tuple of points $(x_{1},\dotsc ,x_{n})$ (with
$x_{i}\not =x_{i+1}$ for all $i\in \Z_{n}=\Z / n\Z$~\footnote{For
notational convenience, we consider the indices modulo $n$.}).  We
can associate a set of weighted configurations~$\C (p)$ on $\partial_\infty
X$ to $p$, by choosing $\xi_{i}$ such that $x_{i+1}\in
\overline{x_{i}\xi_{i}}$, and setting $m_{i}:=d(x_{i},x_{i+1})$.
Then all $c\in \C (p)$ are semistable by 
\cite[Lemma 4.3]{kapovichLeebMillson} (their proof generalizes without
problems). 
Observe that (if $X$ is not geodesically complete) it may happen that
$\C (p)=\emptyset$.

An element $c\in \C (p)$ is
called a \emph{Gauss map} for $p$ (since this construction, in the
case of the hyperbolic plane, was
mentioned in a letter from Gauss to Bolyai, \cite{gauss}).

On the other hand, consider a weighted configuration $c$. 

Let 
\[
\Phi_{c}:=\phi_{\xi_{n},m_{n}}\circ \dotsb \circ \phi_{\xi_{1},m_{1}}.
\]
Since a composition of 1-Lipschitz maps is 1-Lipschitz, $\Phi_{c}$ is
1-Lipschitz, i.e.\ a weak contraction. Every fixed point of $\Phi_{c}$
is a first vertex of a polygon $p$ with $c\in \C (p)$.

\medbreak

A more general
discussion of measures on $\partial_\infty X$ (if $X$ is a symmetric
space or \EB ) can be found in
\cite{KLM1}, \cite{kapovichLeebMillson}.

\begin{figure}[!b]
\hfill 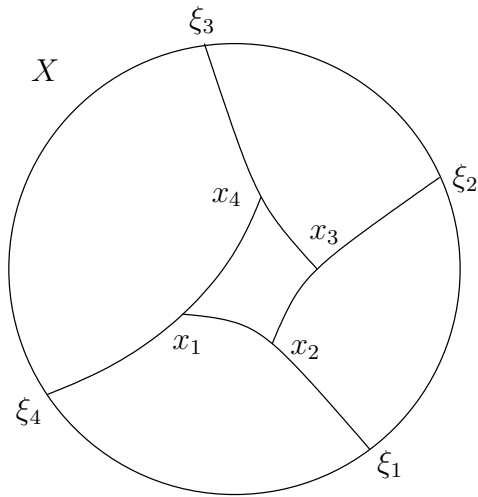 \hfill ~
\caption{Gauss maps}
\end{figure}

\section{Ultralimits, ultraproducts, and asymptotic tubes}
\subsection{Ultralimits}\label{ultralimitsSection}
This section introduces the notion of ultralimit, and the special
cases ultraproduct and asymptotic tube, which play
an important role in our proof.

We keep the general discussion of ultralimits brief and refer the
interested reader to
  \cite[pp.\ 77-80]{bridsonHaefliger} and
\cite[sect.\ 2.4]{kleinerLeeb} for more details.

\begin{definition} Let $\omega$ be a (fixed) non-principal
ultrafilter\footnote{In our context, a non-principal ultrafilter is a
means of (consistently) choosing an accumulation point for any bounded
sequence of real numbers.}, and let $(X_{i},d_{i},o_{i})_{i}$ be a
sequence of metric spaces with metrics $d_{i}$ and basepoints $o_{i}$.  

Then $X_{\omega}:=\lim_{\omega}(X_{i},d_{i},o_{i})$ is the \emph{ultralimit} of
this sequence, a space consisting of
equivalence classes of sequences $(x_{i})$ with $x_{i}\in X_{i}$ and
$d(x_{i},o_{i})$ bounded. 
The distance between two such
sequences $(x_{i,n})_{n}$ (for $i\in \{1,2 \}$) is
$\lim_{\omega}d(x_{1,n},x_{2,n})$, the accumulation point of
$(d(x_{1,n},x_{2,n}))_{n}$ picked by~$\omega$. The equivalence classes
consist of sequences having distance zero.
\end{definition}

If all $X_{i}$ are CAT(0), then their ultralimit is a Hadamard
space; if all $X_{i}$ are (additionally) geodesically complete, then
every geodesic segment, ray and line in $X_{\omega}$ arises as
ultralimit of geodesic segments, rays, and lines respectively
(\cite[2.4.2, 2.4.4]{kleinerLeeb}).

If all $X_{i}$ are \EB s with isometric spherical Weyl chamber, then 
their ultralimit  is also a \EB \
with the same spherical Weyl chamber (\cite[sect. 5.1]{kleinerLeeb}). 

\bigbreak

Let us assume for the rest of this section that
$(X_{i},d_{i})_{i}=(X,d)_{i}$ is a constant 
sequence, and  $X$ is a Hadamard space; so only the basepoint varies in the
construction of the ultralimit $X_{\omega}$. 

Then there is a natural map $*:\partial_\infty X\rightarrow \partial_\infty
X_{\omega}$, obtained by assigning to $\xi \in \partial_\infty X$ the
equivalence class of rays in $X_{\omega}$ which has finite distance
from the ray defined by the sequence of rays $\overline{o_{i}\xi}$. We
denote the image of $\xi$ by $\xi_{*}$. 

Now we can push a weighted configuration $c$ on $\partial_\infty X$ forward to a 
weighted configuration $c_{*}$ on $\partial_\infty X_{\omega}$ by mapping the
$\xi_{i}$ to $\xi_{i,*}$ and keeping the weights. 

\begin{lemma}\label{contractionOnUltralimits}
Under the assumptions above, let $\Phi_{*}$ denote the weak
contraction associated to the pushed forward configuration. Then
$\Phi_{*}$ has the form
\[
\Phi_{*}\Kr{(x_{i})_{i}} = (\Phi (x_{i}))_{i}
\]
\end{lemma}
\begin{proof}
It suffices to show that for any $\xi \in \partial_\infty X$ and a
real number $m>0$, pushing towards $\xi_{*}$ by $\phi_{\xi_{*},m}$ has
the form given above. So let $x=(x_{i})_{i}\in X_{\omega}$. Recall
that by definition,
the distances $d(x_{i},o_{i})$ are bounded. Hence, the ray
$\overline{x\xi_{*}}$ can be represented by the ultralimit of the rays
$\overline{x_{i}\xi}$, which implies the claim.
\end{proof}

\subsection{Ultraproducts}

\begin{definition}
For a metric space $X$ let the
\emph{ultraproduct} of $X$ be the ultralimit of the
constant sequence $(X_{i},d_{i},o_{i}):=(X,d,o)$; i.e.\ 
$X^{\omega}:=\lim_{\omega}(X,d,o)$
 (where we have chosen
a basepoint $o$ for $X$, which has no influence on the isometry type
of $X^{\omega}$).
\end{definition}

There is a canonic isometric embedding $X\rightarrow X^{\omega}$ sending $x$ to
$(x,x,\dotsc )$.
    
Observe that if $X$ is proper (e.g.\ a locally compact CAT(0)-space), the
ultraproduct~$X^{\omega}$ 
is isometric to $X$. 

For details on ultraproducts, see
\cite[sect.\ 11]{lytchakSphB}.

\subsection{Asymptotic tubes}\label{asymptoticTube}

One of the main ideas in the proof of our main theorem is that the
weak contraction $\Phi_{c}$ associated to a weighted configuration
asymptotically moves a ray to a parallel ray.

 We make this idea precise by using particular
ultra-limits.

Throughout this section,  $X$ will be a Hadamard space and $\rho
=\overline{o\eta}$ will be a ray in $X$.

Let $\xi \in \partial_\infty X$.
The following lemma says that pushing towards $\eta$ and $\xi$
asymptotically commutes when moving out along $\rho$.

\begin{lemma}\label{commutativity}
Let $m,c> 0$ and $\xi \in \partial_\infty X$.
Then $\displaystyle \lim_{t\rightarrow \infty}d(\phi_{\xi,m}\circ 
\phi_{\eta ,c}\circ \rho (t), 
\phi_{\eta ,c}\circ \phi_{\xi ,m}\circ \rho (t))=0$.
\end{lemma}
\begin{proof}
Let $o_{t}:=\rho (t), x_{t}:=\phi_{\eta ,c}(o_{t}) = \rho (t+c),
y_{t}:=\phi_{\xi ,m}(o_{t})$, and
 $\hat{\alpha}:=\angle(\eta,\xi )$. We may assume 
$\hat{\alpha}\not =0$, since otherwise $\eta =\xi$, and there is
nothing to show. 

If we set $z_{t}:=\phi_{\eta ,c}(y_{t})$, then the claim is
$d(z_{t},y_{t+c})\underset{t\rightarrow \infty}{\rightarrow} 0$. 

Let $\varepsilon >0$ be given.
\begin{figure}[!b]
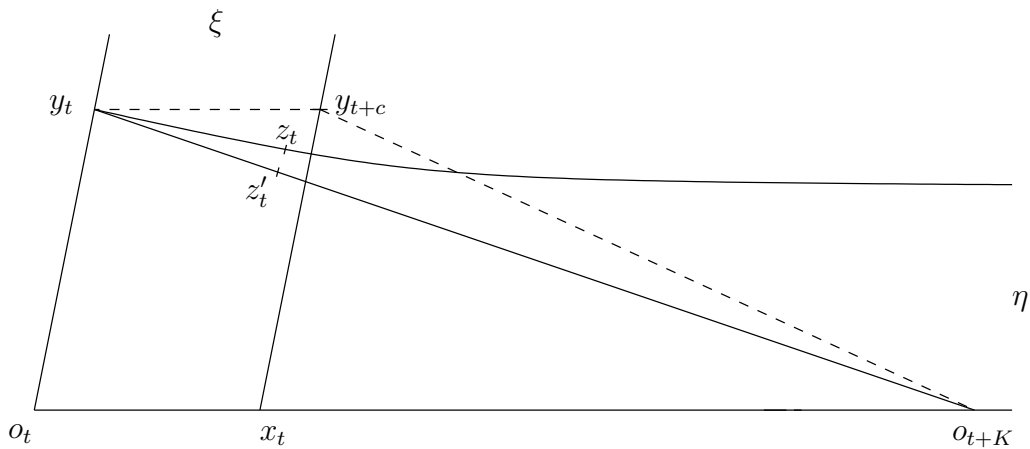
\caption{The points from the proof of Lemma \ref{commutativity}}
\end{figure}
\begin{enumerate}
\item Let $K=K(\varepsilon ,m,c)$ be the constant from Lemma
\ref{rayLemma}. We may assume $K\geq c$.

Let $z'_{t}:=\overline{y_{t}\rho(t+K)}(c)$.
We have $d(z'_{t},z_{t})\leq \varepsilon$, so we try to
get information about $d(z'_{t},y_{t+c})$.

\item 
Let $\bar{\alpha}<\hat{\alpha}$ be such that for a Euclidean triangle
$ABC$ with sides $AC,AB$ of length $K-c,m$ respectively, the length of
the third side varies by at most $\varepsilon$ when the angle at $A$
varies in the interval $[\bar{\alpha},\hat{\alpha}]$.

 Let $l$ be the
maximal length of the third side (occurring when the angle is equal
to~$\hat{\alpha})$.

\item 
Observe that in a Euclidean triangle $ABC$ with sides $AC,AB$ of
length $K,m$ respectively, and angle at $A$ in the interval
$[\bar{\alpha},\hat{\alpha}]$, the third side has length at least
$c+(l-2\varepsilon)$. 

Since the constant $K$ from Lemma
\ref{rayLemma} is independent from the \HS \ (so we may choose
$X=\R^{2}$ here), the claim follows from (2).

\item 
 Finally, let $T>0$ be such that for $t>T$, we have
$\bar{\alpha}\leq \alpha_{t} :=\angle_{o_{t}}(\eta ,\xi )\leq \hat{\alpha}$
(observe that the second inequality is trivial). 

\end{enumerate}
Now we consider the triangle $\Delta (y_{t}y_{t+c}\rho (t+K))$ for
$t>T$. \\
Since the angle corresponding to $\alpha$ in the comparison triangle
is in the interval $[\bar{\alpha},\hat{\alpha}]$, (2)~implies
$d(y_{t+c},\rho (t+K))\leq l$; and since $\phi_{\xi ,m}$ is 
1-Lipschitz, we have $d(y_{t},y_{t+c})\leq c$. On the other hand, we
have $d(y_{t},\rho (t+K))\geq c+l-2\varepsilon$ by (3).

Considering the Euclidean comparison triangle, this shows that we have
control over $d(z'_{t},y_{t+c})$, and this quantity becomes arbitrarily
small as $\varepsilon$ goes to 
zero. With (1), this finishes the proof. 
\end{proof}
\begin{definition}
In the situation described above, define the
CAT(0)-space 
\[
X_{\omega}:=\lim\nolimits_{\omega}(X,d,\rho (i)).
\]
Observe that in $X_{\omega}$, the image of $\rho$ is a line $l$.
Let $T_{\eta} = T_{\rho} :=P_{l}$, and call this space the \emph{asymptotic tube}
of $\eta$ (it is easy to see that $T_{\rho}\cong T_{\rho '}$ if $\rho$
and $\rho '$ are rays asymptotic to $\eta$).
\end{definition}

Consider the map $*:\partial_\infty X\rightarrow \partial_\infty X_{\omega}$ introduced at
the end of section 
\ref{ultralimitsSection}.

\begin{lemma}
We have $*:\partial_\infty X\rightarrow \partial_\infty T_{\eta}$, \newline
and
for any $\xi \in 
\partial_\infty X$, we have $\angle (\xi ,\eta )=\angle (\xi_{*},\eta_{*})$.
\end{lemma}
\begin{proof}
Let $\xi \in \partial_\infty X$ and $m>0$. We claim that the map 
\[
l_{\xi ,m}: 
t\mapsto  (\phi_{\xi ,m}\circ \rho (i+t))_{i} 
\]
defines a line
parallel to $l$ in $X_{\omega}$ (for given $t$, we set the coordinates with
$i+t<0$ arbitrarily; since these are finitely many, they have no
influence on the point defined in~$X_{\omega}$): Indeed, by the Lemma above, the
following equality holds  in $X_{\omega}$ (for $t'>t$):
\[
l_{\xi ,m}(t')=(\phi_{\xi ,m}\circ
\phi_{\eta ,t'-t}\circ \rho (i+t))_{i} = (\phi_{\eta ,t'-t}\circ 
\underbrace{\phi_{\xi ,m}
\circ \rho (i+t)}_{\text{defining }l_{\xi ,m}(t)})_{i}.
\]
The right hand side shows $d(l_{\xi ,m}(t'),l_{\xi ,m}(t))=t'-t$
 for $t'\geq t$; hence, $l_{\xi ,m}$ is a
geodesic line.
Clearly,  $l_{\xi ,m}$ stays within bounded distance of $l$, so it is
parallel to $l$ (by \cite[II.2.13]{bridsonHaefliger}).

For given $t$, we have $\angle_{\rho (i+t)}(\xi ,\eta
)\underset{i\rightarrow \infty}{\rightarrow} \angle (\xi ,\eta )$ 
and $\angle_{\rho (i+t)}(\rho (0),\xi )\rightarrow \pi -\angle (\xi
,\eta )$ (by \cite[Prop.\ 4.2]{ballmann}), so
we find
 $d(l_{\xi ,m},l)=m\sin\angle (\xi ,\eta )$.

It is clear that the flat strip spanned by $l_{\xi ,m'}$ and $l$
contains $l_{\xi ,m}$ for $m'>m>0$, so $\xi$ determines a half-plane
in $P_{l}$ if $\angle (\eta ,\xi )\not =0,\pi$. In the other cases,
$l=l_{\xi ,m}$. 
\end{proof}
\bigbreak

The following  observation is an immediate consequence of the previous
lemma:

\begin{lemma}\label{foldedSlope}
Let $c$ be a weighted configuration on $\partial_\infty X$, and consider the
map $*:\partial_\infty X\rightarrow \partial_\infty T_{\rho}$. Then
$\slope_{c}(\eta )=\slope_{c_{*}}(\eta_{*})$. \qed 
\end{lemma}

\begin{remark}\label{asymptoticTubeRemark}
One can show that $*$ also has the following properties:

 The half-planes determined by  $\xi ,\xi '$ agree if the geodesic
segments $\overline{\eta \xi},\overline{\eta \xi '}$ start in the same
direction. 
The induced map between the spaces of directions
$S_{\eta}(\partial_\infty X)\rightarrow  
S_{\eta_{*}}(\partial_\infty T)$ is 1-Lipschitz, but not an isometric
embedding in general.

We show below that in a \EB , one even gets a map (with the properties
we need)
$*:\partial_{\infty}X\rightarrow \partial_{\infty}P_{l}$ for a line
$l$ containing $\rho$. The same result holds for symmetric spaces of
noncompact type. 

The question arises whether in a general \HS , one can get a suitable map to
the boundary of
$\R \times X_{\eta}$, the space of parametrized strong asymptote
classes at $\eta$ (see \cite{karpelevic}, \cite[sect.\ 2.1.3]{leebHabil},
\cite[sect.\ 3.1.2]{KLM1}). 

However, consider the following subset of the Euclidean plane: 
\[
X=\{(x,y)\st x\geq 1, y\geq \log x \}
\]
With the induced length metric, $X$ becomes a Hadamard space; the
boundary at infinity is an arc of length $\frac{\pi}{2}$. Consider the
boundary point $\eta$ corresponding to the ray $\rho$ in $X$ which is given
by parametrizing the graph of the logarithm with unit speed. Then
$X_{\eta}$ consists of one point only (every ray asymptotic
to $\eta$ eventually lies on the graph of the logarithm), but
$T_{\eta}$ is a half-plane.
\end{remark}

\subsection{Asymptotic tubes in \EB s}
In the case where $X$ is a \EB \ or a symmetric space, the
construction described above specializes to the folding map described
in \cite[sect.\ 3.2.5]{KLM1}.  
We discuss the building case:

\begin{lemma}\label{buildingTube}
Let $X$ be a \EB , $\rho = \overline{o\eta}$ a ray in $X$, and $l$ a
line extending 
$\rho$. Let~$T$ be the asymptotic tube associated to $\rho$. Then
there is a natural isometric embedding $\iota :P_{l}\rightarrow T$,
and we have $Im(*)\subset \partial_{\infty}(\iota (P_{l}))$.
\end{lemma}
\begin{proof}
We state an explicit formula for $\iota$: 
We map $p\in P_{l}$ 
 to $(\phi_{\eta ,i}(p))_{i}$.

Since $\phi_{\eta,t}|_{P_{l}}$ is an isometry of $P_{l}$ for every
$t\geq 0$, the first claim holds.

Let $\xi \in \partial_\infty X$ be a boundary point of $X$. For $t$ large
enough, the rays $\overline{\rho (t)\eta }$ and
$\overline{\rho (t)\xi}$ bound a Euclidean sector (by discreteness of
the angle, see \cite[Axiom 4.1.2.EB2]{kleinerLeeb}). This shows that
$\phi_{\xi ,m}$ eventually maps the ray $\rho$ to a parallel ray.
Since this ray eventually coincides with a line parallel to $l$ by
Lemma \ref{parallels},  the claim follows.
\end{proof}

An immediate consequence is:
\begin{lemma}\label{foldingAndAction}
Let $c$ be a weighted configuration on the boundary of the
\EB~$X$. Let $l$ be a line with 
$\lim_{t\rightarrow \pm \infty}l(t)=\xi_{\pm}$. Let $c_{*}$ denote the weighted
configuration on $\partial_\infty P_{l}$ obtained from $c$ via Lemma
\ref{buildingTube}. Then there exists $T>0$ such that 
\[
\forall t>T: \Phi_{c}\circ l(t) = \Phi_{c_{*}}\circ l(t).
\]
\end{lemma}
\begin{proof}
In the proof of Lemma \ref{buildingTube}, we showed that the
definition of
$\iota$ implies  that the
claim holds for 
configurations consisting of a single point, i.e.\ for maps
$\phi_{\xi,m}$.  

Since $\Phi_{c},\Phi_{c_{*}}$ are finite compositions of such maps, the
lemma follows.
\end{proof}

For \EB s, we obtain the following refinement of Lemma \ref{foldedSlope}:

\begin{lemma}\label{slopeEqual}
Let $X$ be a \EB , and let $c$ be a weighted configuration on its
boundary at infinity. Let $\eta \in \partial_{\infty}X$, and $l$ a line asymptotic
to $\eta$. Consider the measure~$c_{*}$ on $\partial_{\infty}P_{l}$
obtained via Lemma \ref{buildingTube}. Then 
\[
\slope_{c}=\slope_{c_{*}}
\]
on a neighborhood of $\eta$.
\end{lemma}
\begin{proof}
Let $U$ be the neighborhood of $\eta$ consisting of points lying in a
common Weyl chamber with $\eta$, and let $\xi\in U, \xi '\in
\partial_{\infty}X$. It follows from the proof of Lemma
\ref{buildingTube} that $\angle (\xi ,\xi ')=\angle (\xi_{*},\xi
'_{*})$, since the triangles $\xi \eta \xi '$ and 
$\xi_{*}\eta_{*}\xi'_{*}$ are isometric (both are spherical, have two
sides of the same length, and have the same angle at $\eta_{(*)}$).
\end{proof}

\begin{lemma}\label{semistableUnderFolding}
Let $X$ be a \EB , and $c$ a semistable configuration on its boundary
at infinity.  Let $\eta \in \partial_{\infty}X$ be a point with
$\slope_{c}(\eta )=0$, and $l$ a line asymptotic 
to $\eta$. Consider the measure $c_{*}$ on $\partial_{\infty}P_{l}$
obtained via Lemma \ref{buildingTube}. Then $c_{*}$ is semistable on
$P_{l}$. 
\end{lemma}
\begin{proof}
The measure $c_{*}$ is supported on the product $l\times CS(l)$, and
\[
\slope_{c_{*}}(\eta_{*})=\slope_{c}(\eta )=0.
\]
Thus for the
antipode $\eta^{-}_{*}$ of $\eta_{*}$, we have
$\slope_{c_{*}}(\eta^{-}_{*})=-\slope_{c_{*}}(\eta_{*})=0$. 

For a point $\xi$ on $\partial_{\infty}P_{l}$ which has distance less
than $\pi$ from $\eta_{*}$, the claim $\slope_{c_{*}}(\xi )\geq 0$
follows from (strict) convexity of the zero-sublevel set of
$\slope_{c_{*}}$ (\cite[Prop.\ 3.1.(ii)]{KLM1}, together with
Lemma \ref{slopeEqual}.
\end{proof}

\section{The Results}

\subsection{Projecting rays to subspaces}

We  examine how rays project to a subspace of a \HS:

\begin{proposition}\label{projectingRays}
Let $X'$  be a Hadamard space and $X\subset X'$ a closed convex
subset. Consider $\eta \in \partial_{\infty}X'$ 
such that $\angle (\eta ,\partial_{\infty}X)<\frac{\pi}{2}$. Let
$o\in X$, $\rho :=\overline{o\eta}$, and $\pi:X'\rightarrow X$ be the
nearest point
projection. Then the segments $\overline{o\,\,(\pi \circ \rho (t))}$
converge to the 
ray $\overline{o\xi}$ (in the cone topology), where $\xi\in
\partial_{\infty}X$ is the unique point with $\angle (\eta ,\xi
)=\angle (\eta,\partial_{\infty}X)$.
\end{proposition}

\begin{proof}
Observe that $\partial_{\infty}X$ is a closed convex subset of
$\partial_{\infty}X'$ (it is even closed in the cone topology); 
since $\angle (\eta ,\partial_{\infty}X)<\frac{\pi}{2}$, the
projection $\xi$ of $\eta$ exists and is unique
(\cite[II.2.6]{bridsonHaefliger}). 

Let $\bar{\alpha}:=\angle (\eta ,\xi )$,
 $c_{t}:=\rho (t)$, $p_{t}:=\pi (c_{t})$, and
$\alpha_{t}:=\tilde{\angle}_{o}(c_{t},p_{t})$.

By considering triangles $D$ of the form $\Delta (o,c_{t},
\overline{o\xi}(t))$, we conclude $d(c_{t},p_{t})\leq
t\sin \bar{\alpha}$ (since the comparison triangle of $D$ has angle at
most $\bar{\alpha}$ at $o$, the CAT(0)-condition gives the upper bound
on $d(c_{t},p_{t})$);
this implies that $\alpha_{t}\leq \bar{\alpha}$ for all $t>0$.

Since $d(c_{t},p_{t})\leq t\sin\bar{\alpha}$, we have $d(o,p_{t})\geq
t(1-\sin\bar{\alpha})$. Thus, for 
$s(1-\sin\bar{\alpha})\geq t$, the same argument as for the
boundedness of $\alpha_{t}$ shows $\alpha_{t}\leq \alpha_{s}$ $(*)$. 

Let $t_{n}:=(1-\sin\bar{\alpha})^{-n}$ for $n\in \N$ (observe that
$\bar{\alpha}\geq \alpha_{t}>0$ as soon as $c_{t}\not \in X$). By what
we have
shown, $\alpha_{t_{n}}$ is an increasing bounded sequence, which
converges to some~$\hat{\alpha}\leq \bar{\alpha}$. 

Given $\varepsilon >0$, let $N$ be such that $\alpha_{t_{N}}\geq
\hat\alpha -\varepsilon$. Then for $t\geq t_{N+1}$ (so $t\in
[t_{n},t_{n+1}]$ for some $n>N$), we have $\hat{\alpha}-\varepsilon \leq
\alpha_{t_{N}}\leq \alpha_{t}\leq \alpha_{t_{n+2}}\leq \hat{\alpha}$ by $(*)$.
Hence $\alpha_{t}\underset{t\rightarrow
\infty}{\rightarrow}\hat{\alpha}$.

We will show 
next that $d(p_{t},\overline{op_{s}})/t\rightarrow 0$ for $s,t$ large;
since
$d(p_{t},o)\geq t(1-\sin\bar{\alpha})$, this implies that the segments
$\overline{op_{t}}$ 
converge to a ray.

For $s(1-\sin\bar{\alpha})\geq t$, let $p_{s,t}$ be the projection of
$c_{t}$ to the 
segment $\overline{op_{s}}$. 
For $\varepsilon >0$, there exists $T$ such that $t\geq T$ implies
$\sin\alpha_{t}\geq  \sin\hat{\alpha}-\varepsilon$. Then for
$s(1-\sin\bar{\alpha})\geq t\geq T$, we have
$d(c_{t},p_{t})\geq t(\sin\alpha_{t})\geq
t(\sin\hat{\alpha}-\varepsilon )$ and  
$d(c_{t},p_{s,t}) \leq t\sin\alpha_{s}\leq t\sin\hat{\alpha}$.

Consider the comparison triangle $\Delta
(c_{t},p_{t},p_{s,t})$. Since $p_{t}$ is the projection of $c_{t}$ to $X$,
its angle at $p_{t}$ is at least $\frac{\pi}2$. 
Hence for the comparison angle
$\gamma_{s,t}:=\tilde{\angle}_{c_{t}}(p_{t},p_{s,t})$, we 
have $\cos \gamma_{s,t}\geq
\frac{\sin\hat{\alpha}-\varepsilon}{\sin\hat{\alpha}}
\underset{\varepsilon
\rightarrow 0}{\rightarrow}1$. 

Thus 
$d(p_{t},p_{s,t})/t\underset{\varepsilon \rightarrow 0, s
(1-\sin\bar{\alpha}) \geq t\geq
T_{\varepsilon}}{\longrightarrow }0$. 
This shows that the segments $\overline{op_{t}}$ converge to a ray
$\overline{o\xi '}$ for some $\xi '\in \partial_{\infty}X$. 

By \cite[L.~2.3.1]{kleinerLeeb}, we have $\angle (\eta ,\xi ')\leq
\liminf_{t\rightarrow \infty}
\tilde{\angle}(c_{t},p_{t})=\hat{\alpha}\leq \bar{\alpha}$. Hence, 
$\hat{\alpha}=\bar{\alpha}$ and $\xi' =\xi$. 
\end{proof}

\begin{proposition}
Let $X'$ be a Hadamard space and $X\subset X'$ a closed convex
subset. Consider $\eta \in \partial_{\infty}X'$, and assume that for
some $o\in X$, the projection of the ray
$\overline{o\eta}$ to $X$ is bounded, i.e.\ there is $ m$ s.t.\ $d(o,\pi
\circ \overline{o\eta}(t))<m$ for all $t>0$. 

Then there exists a point $p\in X$ s.t.\ $\pi \circ
\overline{p\eta}(t)=p$ for all $t>0$.
\end{proposition}

\begin{proof}
Let $c_{t}:=\overline{o\eta}(t)$ and $p_{t}:=\pi(c_{t})$.

Let $t_{1}:=1$, and define $t_{n}$ inductively by
$t_{n}:=K(\frac{1}{n},m,t_{n-1})$, where $K$ is the constant from Lemma
\ref{rayLemma}. Observe that $t_{n}$ is strictly increasing and unbounded.

Observe that $\pi
(\overline{p_{t_{n}}c_{t_{n}}}(t_{n-1}))=p_{t_{n}}$. Since $\pi $ is
1-Lipschitz, we get from Lemma \ref{rayLemma} that $d(p_{t_{n}},
\pi (\overline{p_{t_{n}}\eta}(t_{n-1})))<\frac{1}{n}$. 

We consider the ultraproducts $X^{\omega}\subset (X')^{\omega}$. 
Let $\pi_{X^{\omega}}:(X')^{\omega}\rightarrow X^{\omega}$ be the
projection. Note that $\pi_{X^{\omega}}$ can be given in the form 
\[
\pi_{X^{\omega}}(x_{n})_{n}=(\pi (x_{n}))_{n}.
\]
Then
$p':=(p_{t_{n}})_{n}$ is a point in $X^{\omega}$ which satisfies
$\pi_{X^{\omega}}(\overline{p'\eta}(t))=p'$ for all $t>0$.

Now let $p$ be the projection of $p'$ to $X$. By the above, we have
$\pi_{X^{\omega}}|_{X'}=\pi$, so $\pi_{X^{\omega}}
(\overline{p\eta}(t))\in X$. On the other hand, $d(\pi_{X^{\omega}}
(\overline{p\eta}(t)), p')\leq d(p,p')=d(p',X)$, so the projection of the
ray $\overline{p\eta}$ is constant.
\end{proof}

\begin{remark}
Observe that a point with the properties from the Lemma above is a
global minimum of the Busemann function $b_{\eta}|_{X}$.

Note also, that the example from Remark \ref{asymptoticTubeRemark}
shows that the assumption 
of the proposition above needs not be fulfilled if $\angle
(\eta,\partial_{\infty}X)\geq \frac{\pi}{2}$.  
\end{remark}

\subsection{Persistence of semistability}

Now persistence of semistability follows easily:

\begin{proposition}
\label{persistence}
Let $X\subset X'$, where $X$ is a closed convex subset of the \HS
~$X'$, and let $c$ be a weighted 
configuration on the asymptotic boundary of $X$. If $c$ is semistable on $X$,
then~$c$ is semistable on $X'$. 
\end{proposition}

\begin{proof}
Assume there is $\eta \in \partial_{\infty}X'$ with
$\slope_{c}(\eta)=-c<0$.  From the formula for the slope, we conclude
that there must be some $\xi_{i}$ in the support of $c$ which
satisfies $\angle (\eta ,\xi_{i})<\frac{\pi}{2}$. Hence
Proposition \ref{projectingRays} applies.

From this point, we obtain a contradiction as in the end of the proof
of \cite[L.~3.10.ii]{KLM1}: 

Use the notation of the proof above, and for $s\geq t$, let
$\bar{p}_{s,t}:= \overline{op_{s}}(\frac{t}{s}d(o,p_{s})$.  
We may normalize $b_{c}$ such
that $b_{c}(o)=0$. Then by convexity, we have $b_{c}(c_{s})\leq -cs$. 
As in the proof of \cite[L.~3.10]{KLM1}, we have $b_{c}\geq b_{c}\circ
\pi $ (where $\pi $ is the projection $X'\rightarrow X$).
In particular, $b(p_{s})\leq -cs$. 

For $s\geq t$, we conclude from convexity that $b_{c}(p_{s,t})\leq -ct$.
Fixing $t$ and letting $s\rightarrow \infty$, this shows
$b_{c}(\overline{o\xi}(t\cos\hat{\alpha}))\leq -ct$, implying
$\slope_{c}(\xi )\leq 
-c/\cos\hat{\alpha}<0$. This is the desired contradiction.
\end{proof}

\begin{remark}
Observe that we cannot expect stability to be preserved under
general embeddings, as one sees e.g.\ by embedding $X$ into $X\times \R$. 

We will only use the above proposition for the inclusion $X\subset
X^{\omega}$. However, we may not expect stability to be preserved in
this case either, as the following example shows:

Consider the disjoint union of copies of $\H^2\times [-n,n]$ for $n\in
\N$, identified along $\H^{2}\times \{0 \}$. This is a \HS ~ by
\cite[II.11.3]{bridsonHaefliger}. Its boundary is precisely the
boundary of $\H^2$, but its ultraproduct contains a copy of $\H^2\times
\R$. 
\end{remark}

\subsection{Proof of the Main Theorem}

In this section we present the proof of our main theorem.

We will need a  lemma about fixed points of weak contractions, which
we recall without proof:

\begin{lemma}[{\cite[Lemma 4.5]{kapovichLeebMillson}}]\label{boundedOrbits}
Let $X$ be a \HS \ of finite diameter. Then every
weak contraction $\Phi :X\rightarrow X$  has a fixed point. \qed 
\end{lemma}

The following lemma was essentially contained in an
earlier version of 
\cite{kapovichLeebMillson}:
\begin{figure}\label{triangleFigure}
\centerline{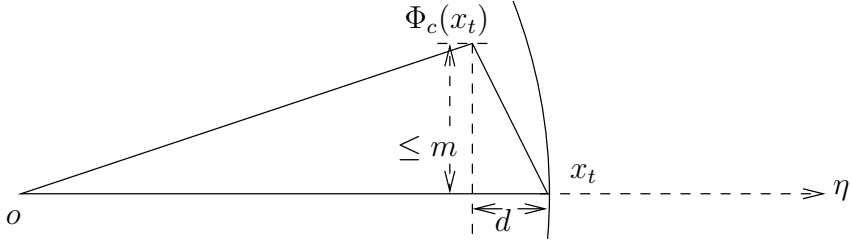}
\caption{For $t$ large, $d(o,\Phi_{c}(x_{t}))<d(o,x_{t})$.}
\end{figure}

\begin{lemma}\label{slopeNonPositive}
Let $c$ be a weighted configuration on the boundary of a \HS\  of the
form $l\times Y$, where $l$ is a line with endpoints $\eta, \eta_{-}$,
and $Y$ is a \HS.

If $\slope_{c}(\eta )>0$, then there exists $T>0$ such that
$d(\Phi_{c}((l(t),y)),(l(0),y))<t$ for all~$t>T$ 
and  $y\in Y$.
\end{lemma}
\begin{proof}
The configuration $c$ can be split into configurations $c_{1},c_{2}$
on $\{\eta ,\eta_{-} \},\partial_\infty Y$ respectively, and this splitting
is compatible with the action of $\Phi$ (see
\cite[L.\ 3.12]{KLM1}). 
In particular, we have $(b_{\eta}\circ \Phi_{c} -
b_{\eta}) \equiv \slope_{c}\eta=:d>0$.

Let $o:=(l(0),y)$ and $x_{t}:=(l(t),y)$.
The triangle $\Delta (o,x_{t},\Phi_{c}(x_{t}))$ is
Euclidean, so the claim follows from the fact that the displacement of
$\Phi_{c}$ is bounded (by $m:=\sum_{i=1}^{n}m_{i}$); see figure
3.
\end{proof}

Now we have all ingredients for the proof of our main theorem; 
we start with the building case:

 \newcounter{lemmaNumber}

\begin{theorem}\label{fixedPointEB}
Let $X$ be a \EB , and let $c$ be a
semistable weighted configuration on its boundary at infinity. Then the
associated weak contraction 
$\Phi_{c}$ has a fixed point.
In particular, there exists a polygon
$p$ in $X$ such that $c$ is a Gauss map for $p$.
\end{theorem}

\begin{proof}
Fix a basepoint $o\in X$.
If we find a ball $B(o,R)\subset X$ which is preserved by $\Phi$, we
are done by Lemma \ref{boundedOrbits}. 

We argue by contradiction: Assume that for each $i\in \N $, there
exists a point $x_{i}\in X$ such that $d(o,x_{i})\geq i$ and 
$ d(\Phi (x_{i}),o) \geq  d(x_{i},o) ~(*)$. Observe that $(*)$ holds
for each $x\in \overline{ox_{i}}$ since $\Phi$ is a weak
contraction. 

The segments $\overline{ox_{i}}$ define a ray $\rho
=\overline{o\eta}$ in the ultraproduct $X^{\omega}$ (for some $\eta
\in \partial_{\infty}X^{\omega}$): We have 
\[
\rho (t)=(\overline{ox_{i}}(t))_{i}
\]
where we set $\overline{ox_{i}}(t):=o$ for $i<t$ (clearly, these
finitely many points have no influence on the point defined in
$X^{\omega}$).

Let $c_{*}$ be the configuration $c$ considered as a configuration on
$\partial_{\infty}X^{\omega}$, and
let $\Phi_{*}$ be the associated weak contraction.
Now $\rho$ satisfies
$d(\Phi_{*}(\rho (t)), o)\geq d(\rho (t),o)=t$ for all $t$, since we have 
\begin{align}
d(\Phi_{*}(\rho (t)), o) = {\textstyle\lim_{\omega}}\underbrace{d(\Phi
(\overline{ox_{i}}(t)),o)}_{\geq t\text{ if }i\geq t}\geq t =d(\rho
(t),o). \tag{$\dagger$}
\end{align}

By Proposition \ref{persistence}, there are two cases to be considered:
\smallbreak

\emph{Case 1: $\slope_{c_{*}}(\eta )>0$:}
We consider the asymptotic tube $T_{\eta}$, and the pushed
forward configuration, which we denote by $c_{**}$; the associated
weak contraction will be denoted by $\Phi_{**}$.

Let $l$ be the line which is obtained from $\rho$ when passing to the
asymptotic tube.
By Lemma \ref{foldedSlope} and Lemma \ref{slopeNonPositive}, we have
$d(\Phi_{**}\circ 
l(t),l(0))<t$ for large $t$. This implies that for large $t$
and $\omega$-almost all $i$, we have 
\(
d(\Phi_{*}\circ \rho (i+t),\rho (i))<t
\).

By the triangle inequality, this implies $d(\Phi_{*}\circ \rho
(i+t),o)<i+t$, in contradiction to~$(\dagger)$.

\smallbreak

\emph{Case 2: $\slope_{c_{*}}(\eta )=0$:}
We argue by induction on rank$(X)$: 

Let $l$ be a line extending $\rho$; we pass to a configuration
$c_{**}$ on $\partial_\infty P_{l}$ (via Lemma \ref{buildingTube}). 
Then $c_{**}$ is semistable by Lemma \ref{semistableUnderFolding}.
Since $P_{l}=l\times CS(l)$, $c_{**}$ splits, and we obtain a
semistable configuration on $\partial_{\infty}l$ and a semistable
configuration on $\partial_{\infty}CS(l)$. 

A semistable configuration on the boundary of a flat Euclidean space
(i.p.\ a line) yields a constant map $\Phi$;
a semistable configuration on $\partial_\infty CS(l)$ has a
fixed point by the induction hypothesis.

Thus, we have a line of fixed points for $c_{**}$ in
$X^{\omega}$. This line of fixed points yields a ray of fixed points
for $\Phi_{*}$ by Lemma \ref{foldingAndAction}. 

So let $p\in X^{\omega}$ be a fixed 
point of~$\Phi_{*}$. There is a unique point $p'\in X$ which is
closest to~$p$. Since $\Phi_{*}$ is 1-Lipschitz, it has to fix $p'$.
Now the observation $\Phi_{*}|_{X}=\Phi$ finishes the proof.
 \end{proof}

\begin{corollary}\label{hadamardCorollary}
Let $X$ be a \HS , and $c$ a
weighted configuration on its boundary at infinity, which  is
stable on $X^{\omega}$. Then the
associated weak contraction 
$\Phi_{c}$ has a fixed point.
In particular, there exists a polygon
$p$ in $X$ such that $c$ is a Gauss map for $p$.
\end{corollary}
\begin{proof}
By assumption, case 2 in the proof of Theorem \ref{fixedPointEB} above
does not occur; hence the proof works exactly the same (observe that
building geometry was used only in the second case).
\end{proof}

In the locally compact case, $X^{\omega}\cong  X$; hence Corollary
\ref{hadamardCorollary}  
 finishes the proof of the Main Theorem.

Observe that we cannot expect Theorem \ref{fixedPointEB} to fully
generalize to \HS s, since in the case of symmetric spaces,  nice
semistability of the configuration is necessary.

\section{Relations to Algebra}

Here, we discuss the relevance of our main theorem to
problems from algebra.
Such problems were studied e.g.\  in \cite{kapovichLeebMillsonAlgebra}.

In the algebraic problems, one only fixes the \emph{type} of a
configuration, i.e.\ the projection  of the points $\xi_{i}$ to the
spherical Weyl chamber $\Delta$. Taking the weights $m_{i}$ into
account, such a type of a configuration may be viewed as an element of
$\Delta^{n}_{euc}$, $n$ copies of the Euclidean Weyl chamber (the
Euclidean cone over the spherical Weyl chamber $\Delta$).
Consider the following theorem:

\begin{theorem}[{\cite[Thm.~1.2]{kapovichLeebMillson}}]
Let $X$ be a \EB . Then for $h\in \Delta^{n}_{euc}$ there exists an
$n$-gon in $X$ with $\Delta$-side lengths $h$ if and only if there
exists a semistable weighted configuration on $\partial_{\infty}X$ of
type $h$.\qed
\end{theorem}

Our main results give a natural proof, and may in 
fact be seen as a refinement, since the proof in \cite{kapovichLeebMillson}
 does not 
provide explicit configurations for which there exists a fixed point.
This indicates that there will eventually be more applications to
algebra.

\bibliographystyle{halpha}

\bibliography{diss}

 \end{document}